\documentclass[12pt]{article}

\usepackage{amssymb}
\usepackage[mathscr]{euscript}

\usepackage{amsfonts}
\usepackage{latexsym}
\usepackage{graphicx}
\newcommand{\bbox}{\normalsize {}%
        \nolinebreak \hfill $\blacksquare$ \medbreak \par}

\author{Fr\'ed\'eric H\'elein}

\newtheorem{theo}{Theorem}
\newtheorem{lemm}{Lemma}

\def\1{\hbox{1\hskip -2pt l}}

\title{Removability of singularities of harmonic maps
into pseudo-Riemannian manifolds}
\author{Fr\'ed\'eric H\'elein}
\large
\begin{document}
\maketitle

{\em Abstract --- We consider harmonic maps into pseudo-Riemannian manifolds. We show the removability of isolated
singularities for continuous maps, i.e.\,that any continuous map from an open subset of $\Bbb{R}^m$ into a
pseudo-Riemannian manifold which is two times continuously differentiable and harmonic everywhere outside an
isolated point is actually smooth harmonic everywhere.}

\section*{Introduction}
Given $n\in \Bbb{N}^*$ and two nonnegative integers $p$ and $q$ such that $p+q=n$,
a pseudo-Riemannian manifold $({\cal N},h)$ of dimension $n$ and of signature $(p,q)$
is a smooth $n$-dimensional manifold ${\cal N}$ equipped with a pseudo-Riemannian metric
$h$, i.e.\,a section of $T^*{\cal N}\odot T^*{\cal N}$ (where $\odot$ is the symmetrised tensor product),
such that $\forall M\in {\cal N}$, $h_M$ is a non degenerate bilinear form of signature $(p,q)$.
Any local chart $\phi:{\cal N}\supset U\longrightarrow V\subset \Bbb{R}^n$ allows us to use local
coordinates $(y^1,\cdots ,y^n)\in V$: we then denote by $h_{ij}(y):=
h_{\phi^{-1}(y)}\left({\partial  \over \partial y^i}, {\partial  \over \partial y^j}\right)$.
We say that $({\cal N},h)$ is of class ${\cal C}^k$ if and only if $h_{ij}$ is ${\cal C}^k$.
We define the Christoffel symbol by
\[
\Gamma^i_{jk}(y):= {1\over 2}h^{il}(y)\left( {\partial h_{lk}\over \partial y^j}(y)
+ {\partial h_{jl}\over \partial y^k}(y) - {\partial h_{jk}\over \partial y^l}(y)\right),
\]
where, as a matrix, $(h^{ij})$ is the inverse of $(h_{ij})$.
Then for any open subset $\Omega$ of $\Bbb{R}^m$ and for any ${\cal C}^2$ map $u$ from $\Omega$ to ${\cal N}$,
if we note
\[
u\simeq \phi\circ u =
\left(\begin{array}{c}u^1\\ \vdots\\u^n\end{array}\right)
\quad \hbox{and}\quad \Gamma_{jk}:= \left(\begin{array}{c}\Gamma^1_{jk}\\ \vdots\\ \Gamma^n_{jk}\end{array}\right)
\]
and if we set $\Delta u:= \sum_{\alpha=1}^m{\partial ^2u\over (\partial x^\alpha)^2}$
and $\Gamma(u)(\nabla u\otimes \nabla u):= \sum_{\alpha=1}^m\Gamma_{jk}(u)
{\partial u^j\over \partial x^\alpha}{\partial u^k\over \partial x^\alpha}$,
we say that $u$ is {\em harmonic into} $({\cal N}, h)$ if and only if
\begin{equation}\label{0.EquaFonda}
\Delta u + \Gamma(u)(\nabla u\otimes \nabla u) = 0.
\end{equation}
Equivalentely we may say that $u$ is a critical point of
\[
{\cal A}[u]:= \int_\Omega h_{ij}(u(x))
\sum_{\alpha=1}^m{\partial u^i\over \partial x^\alpha}{\partial u^j\over \partial x^\alpha}dx^1\cdots dx^m.
\]
The purpose of this paper is to prove the following result. 
\begin{theo}\label{0.MainTheo}
Let $({\cal N},h)$ be a pseudo-Riemannian manifold of class ${\cal C}^2$, $\Omega$ an open subset of
$\Bbb{R}^m$, where $m\geq 2$, $a\in \Omega$ and $u$ a map from $\Omega$ to ${\cal N}$ such that
\begin{itemize}
\item $u$ is continuous
\item $u$ is ${\cal C}^2$ and harmonic on $\Omega\setminus \{a\}$
\end{itemize}
Then $u$ is ${\cal C}^2$ and harmonic on $\Omega$.
\end{theo}
Such a result would be a consequence of standard results if the map $u$ had a finite energy and if $({\cal N},h)$
was Riemannian: indeed one could prove then that $u$ is weakly harmonic (because the capacity of a point vanishes)
and obtain the same conclusion by using the continuity of $u$, thanks to results in \cite{Morrey} and \cite{LU}
(with present form due to S. Hildebrandt).
In dimension 2 the same {\em finite energy} and {\em Riemannian target} hypotheses lead to the same conclusion
but without using the fact that $u$ is continuous as proved in \cite{SacksUhlenbeck}.  However the difficulty here
comes from the fact that the target manifold is pseudo-Riemannian. In particular even if we would assume that the
map $u$ had a finite energy, it would not help much.\\

\noindent
This result answers a question
posed by F. Pedit. It is related to the construction of spectral curves associated to any torus in the
sphere $S^4$, a work in progress by F. Burstall, D. Ferus, K. Leschke, F. Pedit and U. Pinkall (see \cite{BFLPP}
for an exposition of these ideas). Using Theorem \ref{0.MainTheo} these authors are able to prove various
results about Willmore surfaces (recall that the right notion of Gauss maps for Willmore surfaces is the
{\em conformal Gauss map} which takes values into a pseudo-Riemannian homogeneous manifold, see e.g.\,\cite{Bryant},
\cite{Helein2} or \cite{BFLPP}).\\

\noindent
{\em Comments on the proof:} our proof is based on applications of the maximum principle.
The strategy consists roughly on the following: on the one hand we construct a smooth harmonic map
which agrees with the initial one on the boundary of a small ball centered at the singularity, on the other hand we
prove a uniqueness theorem for harmonic maps which takes values in a neighbourhood of a point. The uniqueness
result follows from the maximum principle Theorem \ref{Theo.4} which is inspired by \cite{JaegerKaul}
(see also \cite{Helein}). This reduces the uniqueness problem to an estimate on solutions of elliptic linear PDE's
on a punctured domain, given in Lemma \ref{5.Lemma}, the result where we exploit the fact that the capacity
of a point vanishes. The existence result is obtained through a fixed point argument in H\"older spaces in
Theorem \ref{3.ExistenceTheo}. However in this result we need a uniform estimate in the H\"older topology
of the inital harmonic map. This is the subject of our key result, Theorem \ref{Theo2}, where the uniform
H\"older continuity is established by using the maximum principle. In the course of this paper we give the
proof of more or less standard results for the convenience of the reader: Lemma \ref{1.LemmaPrelim} is
classical in Riemannian geometry, Lemmas \ref{Lemma3} and \ref{5.Lemma} are certainly
well-known to specialists but I did not find proofs of them in the litterature.

\section{Preliminaries}
\subsection{An adaptated coordinates system}
\begin{lemm}\label{1.LemmaPrelim}
Let $({\cal N},h)$ be a pseudo-Riemannian manifold of class ${\cal C}^2$. Let $M_0$ be a point in ${\cal N}$,
$U$ be an open subset of ${\cal N}$ which contains $M_0$ and $\psi:U\longrightarrow W\subset \Bbb{R}^n$
a local ${\cal C}^2$ chart.
Then there exists an open neighbourhood $U_0\subset U$ of $M_0$ with the following property (denoting by $W_0:=\psi(U_0)$):
$\forall M\in U_0$ there exists a smooth diffeomorphism $\Phi_M:W_0\longrightarrow V_M\subset \Bbb{R}^n$
such that in the local coordinates $y \simeq \phi:= \Phi_M\circ \psi\in V_M$,
\begin{equation}\label{1.y=0}
\phi(M) = 0,
\end{equation}
\begin{equation}\label{1.h=eta}
h_{ij}(0) = \eta_{ij} := \hbox{diag}(1,\cdots ,1,-1,\cdots ,-1),
\end{equation}
\begin{equation}\label{1.gamma=0}
\Gamma^i_{jk}(0) = 0,
\end{equation}
and there exists constants $C_W,C_\Gamma>0$ independent of $M\in U_0$ such that
\begin{equation}\label{1.uniformly}
\forall M\in U_0, \quad ||d\Phi_M||_{L^\infty(W_0)} +  ||d\Phi_M^{-1}||_{L^\infty(V_M)}\leq C_W
\end{equation}
and
\begin{equation}\label{1.dgamma}
||d\Gamma||_{L^\infty(V_M)}\leq C_\Gamma.
\end{equation}
\end{lemm}
{\em Proof} --- We denote by $z=(z^1,\cdots ,z^n)\in W$ the coordinates in the local chart 
$\psi:U\longrightarrow W$ and by
$\widetilde{h}_{ij}(z):= h_{\psi^{-1}(z)}\left({\partial  \over \partial z^i},
{\partial  \over \partial z^j}\right)$ the expression of the metric. We look at
a neighbourhood $W_0$ of $\psi(M_0)$ such that for all $M\in \psi^{-1}(W_0)$ there exists 
a map $\Phi_M:W_0\longrightarrow \Bbb{R}^n$ such that, denoting by $z_M:=\psi(M)$,
\[
\forall y\in \Bbb{R}^n,\forall z\in W_0,\ y = \Phi_M(z)\quad \Longleftrightarrow \quad 
z^i - z_M^i = A^i_jy^j + {1\over 2}B^i_{jk}y^jy^k,
\]
for some invertible matrix $A=(A^i_j)\in GL(n,\Bbb{R})$ and real coefficients $B^i_{jk}$ satisfying
$B^i_{jk} = B^i_{kj}$.
This map is well-defined if we choose $W_0$ to be a sufficiently small neighbourhood of $z_M$.
Let us compute the expression $h_{ij}(y)$ of the metric 
in the coordinates $y^i\simeq T^i_M\circ \psi$ in terms of $\widetilde{h}_{ij}(z)$:
\[
h_{kl}(y) = \widetilde{h}_{ij}(z(y)){\partial z^i\over\partial y^k}{\partial z^j\over \partial y^l}
 = \widetilde{h}_{ij}(z(y))\left( A^i_k + B^i_{kp}y^p\right) \left( A^j_l + B^j_{lq}y^q\right).
\]
In order to achieve (\ref{1.h=eta}) it suffices to choose $A^j_j$ such that
$\widetilde{h}_{ij}(z_M)A^i_kA^j_l = \eta_{kl}$ (requiring also that $A$ is symmetric
and positive definite ensures uniqueness). Next we compute that:
\[
{\partial h_{kl}\over \partial y^p}(0) =
{\partial \widetilde{h}_{ij}\over \partial z^r}(z_M)A^r_pA^i_kA^j_l + 
\widetilde{h}_{ij}(z_M)\left( B^i_{kp}A^j_l + A^i_kB^j_{lp}\right),
\quad \forall k,l,p.
\]
And we deduce $\Gamma^i_{jk}$ in function of
$\widetilde{\Gamma}^i_{jk}:= {1\over 2}\widetilde{h}^{il}\left( {\partial \widetilde{h}_{lk}\over \partial z^j}
+ {\partial \widetilde{h}_{jl}\over \partial z^k} - {\partial \widetilde{h}_{jk}\over \partial z^l}\right)$
at 0:
\[
\Gamma^i_{jk}(0) = \eta^{il}\left(\widetilde{h}_{pq}(z_M)\widetilde{\Gamma}^p_{rs}(z_M)A^q_lA^r_jA^s_k
+ \widetilde{h}_{pq}(z_M)A^p_lB^q_{jk}\right).
\]
Since $A^p_l$ and $\widetilde{h}_{pq}(z_M)$ are of rank $m$ and thanks to the
relation $\widetilde{\Gamma}^p_{rs} = \widetilde{\Gamma}^p_{sr}$ we deduce that there exist
unique coefficients $B^i_{jk}$ satisfying $B^i_{jk} = B^i_{kj}$ such that (\ref{1.gamma=0}) holds.
Since $\Phi_M$ depends analytically on $ \widetilde{h}_{ij}(z_M)$
and $\widetilde{\Gamma}^i_{jk}(z_M)$, Condition (\ref{1.uniformly}) is obtained by choosing $W_0$
sufficiently small. Then (\ref{1.dgamma}) is a consequence  of (\ref{1.uniformly}).
\bbox

\subsection{Notations}
In the next two sections we will use spaces of H\"older continuous functions and of
functions with higher derivatives which are H\"older continuous.
We first recall some notations and results from \cite{GilbargTrudinger}.
For any point $x\in \Bbb{R}^m$ or $y\in \Bbb{R}^n$ we note
$|x|:= \sqrt{\sum_{\mu=1}^m(x^\mu)^2}$ and $|y|:= \sqrt{\sum_{i=1}^n(y^i)^2}$.\\

\noindent
For any $\alpha\in (0,1)$ and for all open subset $\Omega\subset \Bbb{R}^m$,
we define ${\cal C}^{0,\alpha}(\Omega)$ to be the set of functions $f$ on $\Omega$ which are
$\alpha$-H\"older continuous on $\Omega$, i.e.\,such that
\[
\hbox{for all compact }K\subset \Omega,\quad \sup_{x,y\in K,x\neq y}
{|f(x) - f(y)|\over |x-y|^\alpha}< \infty.
\]
For $k\in \Bbb{N}$, we let ${\cal C}^{k,\alpha}(\Omega)$ to be the set of ${\cal C}^k$ functions
$f$ on $\Omega$ such that $D^k f\in {\cal C}^{0,\alpha}(\Omega)$. Here
$\forall \mu=(\mu_1,\cdots ,\mu_m)\in \Bbb{N}^m$ we write $|\mu|:= \mu_1+\cdots +\mu_m$ and 
\[
D^kf := \left({\partial ^kf\over (\partial x^1)^{\mu_1} \cdots (\partial x^m)^{\mu_m}}\right)
_{\mu_1,\cdots ,\mu_m\in \Bbb{N}, |\mu| = k}.
\]
We also will use the notation
\[
|D^kf| := \sum_{\mu_1,\cdots ,\mu_m\in \Bbb{N}, |\mu| = k}
\left| {\partial ^kf\over (\partial x^1)^{\mu_1} \cdots (\partial x^m)^{\mu_m}}\right|.
\]
For $x,y\in \Omega$, we denote by $d_x:=\hbox{dist}(x,\partial \Omega)$ and $d_{x,y}:= \hbox{min}(d_x,d_y)$ and
if $\alpha\in (0,1)$, $\beta\in \Bbb{R}$, $k\in \Bbb{N}$, we set
\[
\begin{array}{l}\displaystyle 
[u]^{(\beta)}_{k,0;\Omega}:= [u]^{(\beta)}_{k;\Omega} := \sup_{x\in \Omega}d_x^{\beta+k}|D^k u(x)|
\\
\displaystyle |u|^{(\beta)}_{k,0;\Omega}:= |u|^{(\beta)}_{k;\Omega}:= \sum_{j=1}^k\ [u]^{(\beta)}_{j;\Omega}
\\
\displaystyle [u]^{(\beta)}_{k,\alpha;\Omega}:= \sup_{x,y\in \Omega}d^{\beta + k+\alpha}_{x,y}
{|D^k u(x) - D^k u(y)|\over |x-y|^\alpha}
\\
\displaystyle |u|^{(\beta)}_{k,\alpha;\Omega}:= |u|^{(\beta)}_{k;\Omega} + [u]^{(\beta)}_{k,\alpha;\Omega}.
\end{array}
\]
We  also define
\[
{\cal C}^{(\beta)}_{k,\alpha;\Omega}:= \{u\in {\cal C}^{k,\alpha}(\Omega)/
|u|^{(\beta)}_{k,\alpha;\Omega}< \infty \}.
\]
\noindent
For all $a\in \Bbb{R}^m$, $r\in (0,\infty)$ we note:
\[
B^m(a,r):=\{x\in \Bbb{R}^m/|x-a|<r\}.
\]
We denote by $\omega_m$ the  measure of the unit ball $B^m(0,1)$.
We will use  the following, which is a consequence of Lemma 6.20, Lemma 6.21 and Theorem 6.22 in
\cite{GilbargTrudinger}: 
\begin{lemm}\label{LemmaGT0}
Let $\alpha,\beta\in (0,1)$ and $\Omega$ be a smooth ${\cal C}^2$ open domain of $\Bbb{R}^m$. Then
there exists a constant $C>0$ which depends only on $n$, $\Omega$ and $\alpha$ such that
for all $f\in {\cal C}^{0,\alpha}(\Omega)$ such that $|f|^{(2-\beta)}_{0,\alpha;\Omega}<\infty$,
there exists a unique map $u\in {\cal C}^{2,\alpha}(\Omega)$ which is solution of
\[
\left\{ \begin{array}{cccl}
\Delta u & = & f &\hbox{on }\Omega\\
u & = & 0 & \hbox{on }\partial \Omega.
\end{array}\right.
\]
Moreover $\exists \delta>0$ such that $u$ satifies the estimate
\[
|u|^{(-\beta)}_{2,\alpha;\Omega} \leq \delta |f|^{(2-\beta)}_{0,\alpha;\Omega}.
\]
\end{lemm}
We also recall that
the interpolation result Lemma 6.35 in \cite{GilbargTrudinger} implies the following: if $j,k\in \Bbb{N}$,
$\alpha,\beta\in (0,1)$, $\gamma\in \Bbb{R}$ then 
\begin{equation}\label{1.interpol}
j+\alpha\leq k+\beta \quad \Longrightarrow \quad
|u|^{(\gamma)}_{j,\alpha;\Omega}\leq \widetilde{C} |u|^{(\gamma)}_{k,\beta;\Omega},
\end{equation}
where $\widetilde{C}$ is a positive constant which depends only on $\alpha,\beta, j,k,n$.

\section{H\"older continuity of the map around $a$}
In this section we prove that any map satisfying the hypotheses of Theorem \ref{0.MainTheo}
is uniformly H\"older continuous in a neighbourhood of $a$.

\begin{theo}\label{Theo2}
Let $({\cal N},h)$ be a pseudo-Riemannian manifold of class ${\cal C}^2$, $\Omega$ an open subset of
$\Bbb{R}^m$, $a$ a point in $\Omega$ and $u$ a continuous map from $\Omega$ to ${\cal N}$ which
is ${\cal C}^2$ and harmonic on $\Omega\setminus \{a\}$. Then there exists $\alpha\in (0,1)$ 
and an open ball $B^m(a,R_1)\subset \Omega$ such that $u$ is in ${\cal C}^{0,\alpha}(B^m(a,R_1))$ and 
$|u|_{0,\alpha;B^m(0,R_1)}$ is bounded.
\end{theo}
{\em Proof} --- We apply Lemma \ref{1.LemmaPrelim} with $M_0 = u(a)$. Let $U_0$ be the open neighbourhood of $u(a)$ and
$\psi:U_0\longrightarrow W_0$ be the local chart as in Lemma \ref{1.LemmaPrelim}. We can
assume w.l.g.\,that $\Omega$ is the ball $B^m(a,2R)$, where, since $u$ is continuous, $R\in (0,\infty)$ can been
chosen in such a way that $u(B^m(a,2R))\subset U_0$.\\

\noindent
For any $x_0\in B^m(a,2R)$, let
$\Phi_{u(x_0)}\circ \psi:U_0\longrightarrow V_{u(x_0)}$ be the local chart centered at $u(x_0)$
given by Lemma \ref{1.LemmaPrelim}. It follows from (\ref{1.gamma=0}) and (\ref{1.dgamma}) that
in the coordinate system given by $\Phi_{u(x_0)}\circ \psi$,
\begin{equation}\label{2.0}
\forall y\in V_{u(x_0)}, \forall \xi\in \Bbb{R}^m,\quad
\left|\Gamma_{jk}(y)\xi^j\xi^k\right|\leq C_\Gamma|y||\xi|^2,
\end{equation}
where $C_\Gamma$ is the same constant as in (\ref{1.dgamma}).
We observe that, by replacing $R$ by a smaller number if necessary and
because of the continuity of $u$ and of (\ref{1.uniformly}), we can suppose that
\begin{equation}\label{2.1}
\forall x,x_0\in B^m(a,2R),\quad
\left|\Phi_{u(x_0)}\circ \psi\circ u(x)\right|\leq \inf \left({1\over 2C_\Gamma},{1\over 4}\right).
\end{equation}
We now fix some $x_0\in B^m(a,R)$ and set
\[
u\simeq \Phi_{u(x_0)}\circ \psi\circ u.
\]
We also let, for $r\in (0,R]$,
\[
||u||_{r,x_0} = ||u||_{L^\infty(B^m(x_0,r))}:= \sup_{x\in B^m(x_0,r)}|u(x)|,
\]
and
\[
\lambda := 2C_\Gamma\sup_{x\in B^m(a,2R)}|u(x)|.
\]
Note that, because of the inclusion $B^m(x_0,R)\subset B^m(a,2R)$ and of (\ref{2.1}),
\begin{equation}\label{2.1.5}
2C_\Gamma||u||_{R,x_0}\leq \lambda \leq 1.
\end{equation}
For any $\nu\in S^{n-1}\subset \Bbb{R}^n$ we consider the following functions on $V_{u(x_0)}\subset \Bbb{R}^n$
(which contains $u(B^m(x_0,R))$):
\[
f_+(y) = \langle \nu,y\rangle + \lambda {|y|^2\over 2}\quad \hbox{ and }\quad
f_-(y) = \langle \nu,y\rangle - \lambda {|y|^2\over 2}.
\]
Using (\ref{0.EquaFonda}) and (\ref{2.0}) we find that on $B^m(x_0,R)\setminus \{a\}$,
\[
\begin{array}{ccl}
-\Delta \left(f_+(u)\right) & = & -\langle \nu,\Delta u\rangle - \lambda \Delta \left({|u|^2\over 2}\right)\\
 & = & \langle \nu,\Gamma(u)(\nabla u\otimes\nabla u)\rangle - \lambda |\nabla u|^2 +
 \lambda \langle u,\Gamma(u)(\nabla u\otimes\nabla u)\rangle \\
 & \leq & - \lambda |\nabla u|^2 + C_\Gamma\left( |\nabla u|^2 +\lambda |u|\,|\nabla u|^2\right)|u|.
\end{array}
\]
Since $||u||_{R,x_0}\leq {1\over 4}$ by (\ref{2.1}) and because of (\ref{2.1.5}) we have $\lambda |u(x)|\leq {1\over 4}$,
$\forall x\in B^m(x_0,R)$. So
\begin{equation}\label{2.2}
-\Delta \left(f_+(u)\right) \leq \left({5\over 4}C_\Gamma||u||_{R,x_0}-\lambda\right)|\nabla u |^2\leq 0,
\quad \hbox{on }B^m(x_0,R)\setminus \{a\}
\end{equation}
Similarly
\begin{equation}\label{2.3}
-\Delta \left(f_-(u)\right) \geq 0, \quad \hbox{on }B^m(x_0,R)\setminus \{a\}
\end{equation}
Now fix $r\in (0,R]$ such that $r\neq |x_0|$ and define $D_\varepsilon$ as follows
\begin{itemize}
\item if $|x_0|<r$, for any $\varepsilon\in [0,r-|x_0|)$, $D_\varepsilon:= B^m(a,\varepsilon)\subset B^m(x_0,r)$
\item if $|x_0|>r$, we assume that $\varepsilon = 0$ and set $D_\varepsilon = D_0:= \emptyset$.
\end{itemize}
And we let $u_+^\varepsilon$ and $u_-^\varepsilon$ be the maps from $B^m(x_0,r)\setminus\overline{D_\varepsilon}$
to $V_{u(x_0)}$ which are the solutions of respectively
\[
\left\{ \begin{array}{cccl}
u_+^\varepsilon & = & ||u||_{R,x_0} + \lambda {||u||_{R,x_0}^2\over 2} & \hbox{ on }\partial D_\varepsilon\\
u_+^\varepsilon & = & f_+(u) & \hbox{ on }\partial B^m(x_0,r)\\
-\Delta u_+^\varepsilon & = & 0 & \hbox{ on }B^m(x_0,r)\setminus\overline{D_\varepsilon},
\end{array} \right.
\]
\[
\left\{ \begin{array}{cccl}
u_-^\varepsilon & = & -||u||_{R,x_0} - \lambda {||u||_{R,x_0}^2\over 2} & \hbox{ on }\partial D_\varepsilon\\
u_-^\varepsilon & = & f_-(u) & \hbox{ on }\partial B^m(x_0,r)\\
-\Delta u_-^\varepsilon & = & 0 & \hbox{ on }B^m(x_0,r)\setminus\overline{D_\varepsilon}.
\end{array} \right.
\]
Since
\[
\left\{ \begin{array}{cccl}
u_+^\varepsilon & \geq & f_+(u) & \hbox{ on }\partial \left( B^m(x_0,r)\setminus\overline{D_\varepsilon}\right)\\
-\Delta u_+^\varepsilon = 0 & \geq & -\Delta f_+(u) & \hbox{ on }B^m(x_0,r)\setminus\overline{D_\varepsilon},
\end{array} \right.
\]
the maximum principle implies that
\[
u_+^\varepsilon \geq f_+(u)  \quad \hbox{ on }B^m(x_0,r)\setminus\overline{D_\varepsilon}.
\]
Now we fix an arbitrary compact $K\subset B^m(x_0,r)\setminus \{a\}$. Then for $\varepsilon$ sufficiently small
we have
\[
u_+^\varepsilon \geq f_+(u)  \quad \hbox{ on }K.
\]
We let $\varepsilon$ goes to 0: since $\{a\}$ has a vanishing capacity, the restriction of $u_+^\varepsilon$ to
$K$ converges in $L^1(K)$ to $u_+:= u_+^0$ (apply Lemma \ref{5.Lemma} to $\phi_\varepsilon:= u_+^\varepsilon - u_+$).
Hence
\[
u_+ \geq f_+(u)  \quad \hbox{ on }K.
\]
Since $u$ and $u_+$ are continuous on $B^m(x_0,r)$ and since $K$ is arbitrary we deduce that
\begin{equation}\label{2.4}
f_+(u) \leq u_+  \quad \hbox{ on }B^m(x_0,r).
\end{equation}
Similarly we get
\begin{equation}\label{2.5}
u_- \leq f_-(u)  \quad \hbox{ on }B^m(x_0,r).
\end{equation}
We deduce from (\ref{2.4}) and (\ref{2.5}) that for any $\nu\in S^m$
\[
\left\{ \begin{array}{rcl}
u_- & \leq & \langle \nu,u\rangle - \lambda {|u|^2\over 2}\\
\langle \nu,u\rangle + \lambda {|u|^2\over 2} & \leq & u_+
\end{array} \right.\hbox{ on }B^m(x_0,r)
\]
and thus
\begin{equation}\label{2.6}
u_- \leq u_- + \lambda {|u|^2\over 2} \leq \langle \nu,u\rangle \leq u_+ - \lambda {|u|^2\over 2}\leq u_+,
\quad \hbox{on }B^m(x_0,r).
\end{equation}
Hence if we let $v:B^m(x_0,r)\longrightarrow V_{u(x_0)}$ be the solution of
\[
\left\{ \begin{array}{cccl}
v & = & u & \hbox{ on }\partial B^m(x_0,r)\\
-\Delta v & = & 0 & \hbox{ on }B^m(x_0,r),
\end{array} \right.
\]
and $Q:B^m(x_0,r)\longrightarrow \Bbb{R}$ be the solution of
\[
\left\{ \begin{array}{cccl}
Q & = & {|u|^2\over 2} & \hbox{ on }\partial B^m(x_0,r)\\
-\Delta Q & = & 0 & \hbox{ on }B^m(x_0,r),
\end{array} \right.
\]
Then, since actually $u_{\pm} = \langle \nu,v\rangle \pm \lambda Q$, (\ref{2.6}) implies that
\[
\forall \nu\in S^{n-1},\quad \left|\langle \nu,u\rangle - \langle \nu,v\rangle \right|
\leq \lambda Q,\quad \hbox{ on }B^m(x_0,r).
\]
Hence since $\nu$ is arbitrary and using the maximum principle for $Q$ we obtain
\begin{equation}\label{2.7}
|u-v|\leq \lambda Q\leq \lambda {||u||^2_{r,x_0}\over 2},\quad \hbox{ on }B^m(x_0,r).
\end{equation}
This implies in particular that, since $u(x_0) = 0$,
\begin{equation}\label{2.8}
|v(x_0)|\leq  \lambda {||u||^2_{r,x_0}\over 2}.
\end{equation}
Moreover since $v$ is harmonic, for all $x\in B^m(x_0,{r\over 2})$ we have
(observing that $B^m(x,{r\over 2})\subset B^m(x_0,r)$)
\[
{\partial v\over \partial x^\mu}(x) = {2^m\over \omega_mr^m}\int_{B^m(x,{r\over 2})}{\partial v\over \partial x^\mu}
= {2^m\over \omega_mr^m}\int_{\partial B^m(x,{r\over 2})}v\nu_\mu ds,
\]
where $\nu$ is the exterior normal vector to the boundary. Hence $\forall x\in B^m(x_0,{r\over 2})$,
\[
\left| {\partial v\over \partial x^\mu}(x)\right| \leq {2\over r}\sup_{\partial B^m(x,{r\over 2})}|v|
\leq {2\over r}\sup_{B^m(x_0,r)}|v| \leq {2\over r}||u||_{r,x_0},
\]
where we used the maximum principle for $v$. Hence we deduce that
\begin{equation}\label{2.9}
\forall x\in B^m(x_0,{r\over 2}),\quad
|v(x)-v(x_0)| \leq 2||u||_{r,x_0}{|x-x_0|\over r}.
\end{equation}
And from (\ref{2.8}) and (\ref{2.9}) we get
\begin{equation}\label{2.10}
\forall x\in B^m(x_0,{r\over 2}),\quad
|v(x)| \leq 2||u||_{r,x_0}{|x-x_0|\over r} + \lambda {||u||_{r,x_0}^2\over 2}.
\end{equation}
Using this inequality together with (\ref{2.7}) we obtain that
\begin{equation}\label{2.11}
\forall x\in B^m(x_0,{r\over 2}),\quad
|u(x)| \leq 2||u||_{r,x_0}{|x-x_0|\over r} + \lambda ||u||_{r,x_0}^2.
\end{equation}
We now choose any $\rho\in (0,{r\over 2}]$ and take the supremum of the left hand side of (\ref{2.11})
over $B^m(x_0,\rho)$. It gives
\begin{equation}\label{2.12}
\forall r\in (0,R],\forall \rho\in (0,{r\over 2}],\quad
||u||_{\rho,x_0} \leq 2||u||_{r,x_0}{\rho\over r} + \lambda ||u||_{r,x_0}^2.
\end{equation}
For $k\in \Bbb{N}$ let $r_k:= R4^{-k}$ and apply (\ref{2.12}) for $r=r_k$ and $\rho = r_{k+1}$:
\[
||u||_{r_{k+1},x_0} \leq {||u||_{r_k,x_0}\over 2} + \lambda ||u||^2_{r_k,x_0}.
\]
This implies, denoting by $a_k:= \lambda ||u||_{r_k,x_0}$, that
\begin{equation}\label{2.13}
a_{k+1} \leq {a_k\over 2} + a_k^2.
\end{equation}
We observe that, because of its definition, $a_k$ is a positive decreasing sequence
and, as a consequence of (\ref{2.1}) and (\ref{2.1.5}),
\begin{equation}\label{2.14}
a_k \leq \lambda ||u||_{R,x_0} \leq {1\over 4}.
\end{equation}
We now let $f:\Bbb{R}\longrightarrow \Bbb{R}$ be the function defined by $f(q) = q^2 - {q\over 2}$
and consider the smooth function $\phi:[0,\infty )\longrightarrow \Bbb{R}$ which is a solution of
\[\left\{
\begin{array}{ccll}
\phi(0) & = & a_0 & \\
\displaystyle {d\phi\over dt} & = & f(\phi) = \phi^2 - {\phi \over 2} & \hbox{on }[0,\infty).
\end{array}\right.
\]

\begin{lemm}\label{Lemma2}
Let $\left( a_k\right)_{k\in\Bbb{N}}$ be a decreasing sequence in $[0,{1\over 4}]$
which satifies (\ref{2.13}) and $\phi$ be defined as above. Then $\forall k\in \Bbb{N}$,
\begin{equation}\label{2.15}
a_k \leq \phi(k).
\end{equation}
\end{lemm}
{\em Proof of the Lemma} --- We show (\ref{2.15}) by induction. This inequality is obviously true for
$k=0$. Let us assume that (\ref{2.15}) is true for some value $k\in \Bbb{N}$. We first observe that,
since $f(0)=0$ and $f\leq 0$ on $[0,{1\over 2}]$, $0\leq a_0\leq {1\over 4}$ implies that
\[
\forall t\in [0,\infty ),\quad 0\leq \phi(t)\leq {1\over 4}.
\]
Hence
\begin{itemize}
\item the fact that $f<0$ on $(0,{1\over 4}]$ implies that $\phi$ is decreasing on $[0,\infty )$
\item the fact that $f$ is decreasing on $[0,{1\over 4}]$ implies that $f\circ \phi$ is increasing
on $[0,\infty )$.
\end{itemize}
Thus
\[
\forall t\in [k,\infty),\quad f\circ \phi(t) \geq f\circ \phi(k).
\]
And hence
\[
\begin{array}{ccl}
\phi(k+1) - \phi(k) & = & \displaystyle \int_k^{k+1}\dot{\phi}(t)\,dt\\
&  = & \displaystyle \int_k^{k+1}f\circ \phi(t)\,dt\\
& \geq & \displaystyle \int_k^{k+1}f\circ \phi(k)\,dt = f\circ \phi(k).
\end{array}
\]
Thus
\[
\phi(k+1) \geq \phi(k) + f(\phi(k)) = {\phi(k)\over 2} + \phi(k)^2.
\]
Now since (\ref{2.15}) is true for $k$, i.e.\,$a_k\leq \phi(k)$, we deduce that
\[
a_{k+1} = {a_k\over 2} + a_k^2 \leq {\phi(k)\over 2} + \phi(k)^2 \leq \phi(k+1).
\]
\bbox
\noindent
{\em Back to the proof of Theorem \ref{Theo2}} --- An easy quadrature shows that
$\phi(t) = {a_0\over 2a_0 + e^{t/ 2}(1-2a_0)}$. Hence Lemma \ref{Lemma2} implies
\[
a_k\leq \phi(k)\leq {a_0\over 1-2a_0}e^{-k/2}
\]
and since $0< a_0 \leq {1\over 4}$, we deduce that $a_k\leq {e^{-k/2}\over 2}$, i.e.
\begin{equation}\label{2.16}
||u||_{R4^{-k},x_0}\leq {1\over 2\lambda} e^{-k/2},\quad \forall k\in \Bbb{N}.
\end{equation}
We now choose any $r\in (0,R]$. Then $\exists !k\in \Bbb{N}$ such that
${R\over 4^{k+1}}< r \leq {R\over 4^k}$. Then on the one hand $r \leq {R\over 4^k}$ implies
\[
||u||_{r,x_0} \leq ||u||_{R4^{-k},x_0} \leq {1\over 2\lambda}e^{-k/2}
\]
by (\ref{2.16}). On the other hand ${R\over 4^{k+1}}< r$ $\Longleftrightarrow$
$k> {\log {R\over r}\over \log 4} -1$ implies
\[
{1\over 2\lambda}e^{-k/2} < {e^{1/2}\over 2\lambda} \left( {r\over R}\right)^{1\over 2\log 4}.
\]
Hence we deduce that
\[
\forall r\in (0,R],\quad ||u||_{r,x_0}\leq {e^{1/2}\over 2\lambda}\left( {r\over R}\right)^{1\over 2\log 4}.
\]
We conclude that, since $\lambda$ is independent of $x_0$, $u$ is uniformly H\"older continuous on
$B^m(a,R)$. \bbox

\section{Existence of a smooth solution around $a$}
We start with the following (classical) preliminary.
\begin{lemm}\label{Lemma3}
Let $\Omega$ be an open subset of $\Bbb{R}^m$ whose boundary is ${\cal C}^2$.
Let $\phi\in {\cal C}^{0,\alpha}(\partial \Omega)$ and $f$ be the solution of
\[
\left\{ \begin{array}{cccl}
-\Delta f & = & 0, & \hbox{on }\Omega\\
f & = & \phi, & \hbox{on }\partial \Omega.
\end{array}\right.
\]
Then $f$ is ${\cal C}^{2,\alpha}$ on $\Omega$ and
\begin{equation}\label{3.00}
[f]^{(-\alpha)}_{0,\alpha;\Omega} + [f]^{(-\alpha)}_{1,\alpha;\Omega} +
[f]^{(-\alpha)}_{2,\alpha;\Omega}\leq C_1|\phi|_{0,\alpha;\partial \Omega},
\end{equation}
where  $C_1$ is a positive constant which depends only on $\Omega$.
\end{lemm}
{\bf Remark} --- We do not have an estimate on $|f|^{(-\alpha)}_{2,\alpha;\Omega}$. Indeed
this quantity is in general infinite because $|f|^{(-\alpha)}_{0,0;\Omega}
=\sup_{x\in \Omega}d_x^{-\alpha}|f(x)|$ cannot be finite unless the trace of $f$ on $\partial \Omega$
vanishes. However the maximum principle and (\ref{3.00}) imply the following:
\begin{equation}\label{3.0.5}
|f|^{(0)}_{0,\alpha;\Omega} \leq \left(1+ C_1\left({\hbox{diam}\Omega\over 2}\right)^\alpha\right)
|\phi|_{0,\alpha;\partial \Omega}.
\end{equation}
{\em Proof} --- {\em First step : $\Omega$ is a half space} --- We assume that
$\Omega = \Bbb{R}^m_+:= \{x = (\vec{x},t)/\vec{x}\in \Bbb{R}^{m-1},t\in (0,\infty )\}$. We
use Proposition 7 and Lemma 4 in Chapter V of \cite{Stein}: there exists a constant $C_0'>0$
such that
\begin{equation}\label{3.01}
\sup_{\vec{x}\in \Bbb{R}^{m-1}}|Df(\vec{x},t)| \leq C_0't^{-1+\alpha}|\phi|_{0,\alpha;\Bbb{R}^{m-1}},
\quad \forall (\vec{x},t)\in \Bbb{R}^m_+.
\end{equation}
Moreover using the fact that $Df$ is harmonic, $\forall x\in \Bbb{R}^m_+$ we have if $\rho:=t/2$,
\[
{\partial Df\over \partial x^\mu}(x) = {1\over \omega_m\rho^m}\int_{B^m(x,\rho)}{\partial Df\over \partial x^\mu}(y)dy
= {1\over \omega_m\rho^m}\int_{\partial B^m(x,\rho)}Df(y)\nu_\mu ds(y),
\]
which implies by (\ref{3.01})
\[
\begin{array}{ccl}\displaystyle 
\left|{\partial Df\over \partial x^\mu}(x)\right| & \leq & \displaystyle 
{1\over \omega_m\rho^m}\int_{\partial B^m(x,\rho)}|Df(y)| ds(y)\\
& \leq & \displaystyle  {m\over \rho}C_0'|\phi|_{0,\alpha;\Bbb{R}^{m-1}}\left({t\over 2}\right)^{-1+\alpha}
= 2^{2-\alpha}mC_0'|\phi|_{0,\alpha;\Bbb{R}^{m-1}}t^{-2+\alpha}.
\end{array}
\]
Hence we obtain that there exists a constant $C_0''$ such that
\begin{equation}\label{3.02}
\sup_{\vec{x}\in \Bbb{R}^{m-1}}|D^2f(\vec{x},t)| \leq C_0''|\phi|_{0,\alpha;\Bbb{R}^{m-1}}t^{-2+\alpha},\quad \forall (\vec{x},t)\in \Bbb{R}^m_+.
\end{equation}
A similar reasonning starting from (\ref{3.02}) leads to
\begin{equation}\label{3.03}
\sup_{\vec{x}\in \Bbb{R}^{m-1}}|D^3f(\vec{x},t)| \leq C_0'''|\phi|_{0,\alpha;\Bbb{R}^{m-1}}
t^{-3+\alpha},\quad \forall (\vec{x},t)\in \Bbb{R}^m_+,
\end{equation}
for some constant $C_0'''>0$.\\
Now using (\ref{3.01}) and (\ref{3.02}) we can estimate $[f]^{(-\alpha)}_{1,\alpha;\Bbb{R}^m_+}$ as follows: if
$x=(\vec{x},t)$ and $y=(\vec{y},s)$ are in $\Bbb{R}^m_+$ let $d:=\inf (t,s)$. Then if
$|x-y|\leq 2d$ we have by (\ref{3.02})
\[
\begin{array}{ccl}\displaystyle 
{|Df(x)-Df(y)|\over |x-y|^\alpha} & \leq & \displaystyle \sup_{\tau>d}|D^2f(\xi,\tau)|\ |x-y|^{1-\alpha}\\
& \leq & \displaystyle  C_0''|\phi|_{0,\alpha;\Bbb{R}^{m-1}}d^{-2+\alpha}(2d)^{1-\alpha}
= 2^{1-\alpha}C_0''|\phi|_{0,\alpha;\Bbb{R}^{m-1}}d^{-1}.
\end{array}
\]
On the other hand if $|x-y|> 2d$ we have by (\ref{3.01})
\[
\begin{array}{ccl}\displaystyle 
{|Df(x)-Df(y)|\over |x-y|^\alpha} & \leq & \displaystyle  {|Df(x)|+|Df(y)|\over d^\alpha}2^{-\alpha}\\
& \leq & \displaystyle  {2C_0'd^{-1+\alpha}|\phi|_{0,\alpha;\Bbb{R}^{m-1}}\over d^\alpha}2^{-\alpha}
= 2^{1-\alpha}C_0'|\phi|_{0,\alpha;\Bbb{R}^{m-1}}d^{-1}.
\end{array}
\]
Thus taking into account both cases we find that $[f]^{(-\alpha)}_{1,\alpha;\Bbb{R}^m_+}\leq
2^{1-\alpha}\sup (C_0',C_0'')|\phi|_{0,\alpha;\Bbb{R}^{m-1}}$. An analog work with (\ref{3.02}) and (\ref{3.03})
instead of (\ref{3.01}) and (\ref{3.02}) leads to $[f]^{(-\alpha)}_{2,\alpha;\Bbb{R}^m_+}\leq
2^{1-\alpha}\sup (C_0'',C_0''')|\phi|_{0,\alpha;\Bbb{R}^{m-1}}$.\\

\noindent
The estimate for $[f]^{(-\alpha)}_{0,\alpha;\Bbb{R}^m_+}$ follows from a slightly different argument.
Again let $x=(\vec{x},t)$ and $y=(\vec{y},s)$ be in $\Bbb{R}^m_+$ and let $d:=\inf (t,s)$. If
$|x-y|\leq 2d$ the same reasoning as above works using (\ref{3.01}) and gives
\begin{equation}\label{3.04}
{|f(x)-f(y)|\over |x-y|^\alpha}  \leq  2^{1-\alpha}C_0'|\phi|_{0,\alpha;\Bbb{R}^{m-1}}.
\end{equation}
However if $|x-y| > 2d$, then we write $|f(\vec{x},t)-f(\vec{y},s)| \leq |f(\vec{x},t)-f(\vec{x},0)| +
|f(\vec{x},0)-f(\vec{y},0)| + |f(\vec{y},0)-f(\vec{y},s)|$ and estimate separately each term.
Using again (\ref{3.01}):
\[
\begin{array}{ccl}\displaystyle 
|f(\vec{x},t)-f(\vec{x},0)| & \leq & \displaystyle \int_0^t |Df(\vec{x},\tau)|d\tau\\
& \leq & \displaystyle \int_0^tC_0'|\phi|_{0,\alpha;\Bbb{R}^{m-1}}\tau^{-1+\alpha}d\tau
= {C_0'\over \alpha}|\phi|_{0,\alpha;\Bbb{R}^{m-1}}t^\alpha.
\end{array}
\]
Similarly one gets
\[
|f(\vec{y},s)-f(\vec{y},0)| \leq {C_0'\over \alpha}|\phi|_{0,\alpha;\Bbb{R}^{m-1}}s^\alpha.
\]
Lastly using $|f(\vec{x},0)-f(\vec{y},0)| \leq |\phi|_{0,\alpha;\Bbb{R}^{m-1}}|\vec{x}-\vec{y}|^\alpha$,
one concludes that
\[
\begin{array}{ccl}\displaystyle 
|f(\vec{x},t)-f(\vec{y},s)| & \leq & \displaystyle |\phi|_{0,\alpha;\Bbb{R}^{m-1}}
\left({C_0'\over \alpha}t^\alpha + {C_0'\over \alpha}s^\alpha + |\vec{x}-\vec{y}|^\alpha\right)\\
& \leq & \displaystyle |\phi|_{0,\alpha;\Bbb{R}^{m-1}}\sup \left(1,{C_0'\over \alpha}\right)
\left(t^\alpha + s^\alpha + |\vec{x}-\vec{y}|^\alpha\right).
\end{array}
\]
Assume for instance that $s<t$, so that $d=s$. Then by the Minkowski inequality
$t^\alpha + s^\alpha = d^\alpha + ((t-s)+d)^\alpha
\leq d^\alpha +(t-s)^\alpha + d^\alpha = 2d^\alpha + (t-s)^\alpha$. Hence
\[
t^\alpha + s^\alpha + |\vec{x}-\vec{y}|^\alpha \leq 2d^\alpha + (t-s)^\alpha + |\vec{x}-\vec{y}|^\alpha
\leq 2d^\alpha + 2|x-y|^\alpha < (2^{1-\alpha} +2)|x-y|^\alpha.
\]
And we thus get
\begin{equation}\label{3.05}
|f(\vec{x},t)-f(\vec{y},s)| \leq (2^{1-\alpha} +2)|\phi|_{0,\alpha;\Bbb{R}^{m-1}}\sup \left(1,{C_0'\over \alpha}\right)|x-y|^\alpha.
\end{equation}
So (\ref{3.04}) and (\ref{3.05}) implies the result on $[f]^{(-\alpha)}_{0,\alpha;\Bbb{R}^m_+}$.\\

\noindent
{\em Step 2 --- estimate on an arbitrary domain} --- If $\Omega$ is a domain with a smooth ${\cal C}^2$
boundary, then using  local chart and a partition of unity one can construct an extension
$g\in {\cal C}^{0,\alpha}(\overline{\Omega})$ of $\phi\in {\cal C}^{0,\alpha}(\partial \Omega)$
which satisfies
\[
[g]^{(-\alpha)}_{0,\alpha;\Omega} + [g]^{(-\alpha)}_{1,\alpha;\Omega} +
[g]^{(-\alpha)}_{2,\alpha;\Omega}\leq C_1'|\phi|_{0,\alpha;\partial \Omega}.
\]
Then the harmonic extension of $\phi$ is $f = g + h$, where $h$ is a function which vanishes on $\partial \Omega$
and which satisfies $-\Delta h = \Delta g$ on $\Omega$. Because of the previous estimate on $g$,
$[\Delta g]^{(2-\alpha)}_{0,\alpha;\Omega}\leq C_1'|\phi|_{0,\alpha;\partial \Omega}$. Now  Lemma \ref{LemmaGT0}
implies that $|h|^{(-\alpha)}_{2,\alpha;\Omega} \leq \delta |\Delta g|^{(2-\beta)}_{0,\alpha;\Omega}$. Hence the estimate
on $f$ follows by summing the estimates on $g$ and $h$. \bbox

\begin{theo}\label{3.ExistenceTheo}
Let $({\cal N},h)$ be a pseudo-Riemannian manifold of class ${\cal C}^2$, $\Omega$ an open subset of
$\Bbb{R}^m$, $a\in \Omega$ and $u$ a continuous map from $\Omega$ to ${\cal N}$ which is ${\cal C}^2$
and harmonic on $\Omega\setminus \{a\}$. Then there exists $\alpha\in (0,1)$ and an open ball
$B^m(a,R_2)$ such that $\forall r\in (0,R_2)$, there exists a map $\underline{u}\in
{\cal C}^{0,\alpha}(\overline{B^m(a,r)},{\cal N})\cap {\cal C}^{2,\alpha}(B^m(a,r),{\cal N})$
which is a solution of
\begin{equation}\label{3.0}
\left\{ \begin{array}{cccl}
\Delta \underline{u} + \Gamma(\underline{u})(\nabla \underline{u}\otimes \nabla \underline{u}) & = &
0 & \hbox{on }B^m(a,r)\\
\underline{u} & = & u & \hbox{on }\partial B^m(a,r).
\end{array}\right.
\end{equation}
I.e.\,$\underline{u}$ is a harmonic map with values in ${\cal N}$ which agrees with $u$ on the boundary of
$B^m(a,r)$. 
\end{theo}
{\em Proof } --- Again we start by applying Lemma \ref{1.LemmaPrelim} with $M_0=u(a)$: it provides
us with a local chart $\Phi_{u(a)}\circ \psi$ on ${\cal N}$ around $u(a)$. We denote by
$y^i$, $h_{ij}$ and $\Gamma^i_{jk}$ respectively the coordinates, the metric and the Christoffel
symbols in this chart. In the following we make the identification
$u\simeq \Phi_{u(a)}\circ \psi\circ u$, so that we view $u$ as a map from $B^m(a,\overline{R})$ to
$\Bbb{R}^n$ such that $u(a) = 0$ and the majorations (\ref{1.gamma=0}) and (\ref{1.dgamma}) hold.\\

\noindent
For every $r\in (0,\overline{R}]$ we let $v:B^m(a,r)\longrightarrow \Bbb{R}^n$ be the harmonic extension
of $u$ inside $B^m(a,r)$, i.e.
\begin{equation}\label{3.1}
\left\{ \begin{array}{cccl}
\Delta v & = & 0 & \hbox{on }B^m(a,r)\\
v & = & u & \hbox{on }\partial B^m(a,r).
\end{array}\right.
\end{equation}
We first apply Theorem \ref{Theo2} which ensures us that the
${\cal C}^{0,\alpha}$ norm of $u$ in a neighbourhood of $u$ is bounded:
$\exists \alpha\in (0,1)$, $\exists R_1\in (0,\overline{R}]$ such that $|u|_{0,\alpha;B^m(a,R_1)}$ is finite.
This allows us to use Lemma \ref{Lemma3} in order to estimate $v$: we will use the notations
\[
|v|^{(-\alpha)}_{k,\alpha;r}:= |v|^{(-\alpha)}_{k,\alpha;B^m(a,r)},\quad \forall k\in \Bbb{N},
\quad \hbox{and}\quad |[v]|_{\alpha;r}:= |v|^{(0)}_{0,\alpha;r} + [v]^{(-\alpha)}_{1,\alpha;r}
+ [v]^{(-\alpha)}_{2,\alpha;r}
\]
and then (\ref{3.00}) and (\ref{3.0.5}) imply
\begin{equation}\label{3.2}
\forall r\in (0,R_1],\quad |[v]|_{\alpha;r}\leq C_2|u|_{0,\alpha;R_1},
\end{equation}
where $C_2 =  C_1+ 1 + C_1\overline{R}^\alpha$. We will denote by
\begin{equation}\label{3.2.5}
\Lambda:= \sup_{r\in (0,R_1]}|[v]|_{\alpha;r}.
\end{equation}
Our purpose is to construct the extension $\underline{u}$ satisfying (\ref{3.0}). By writing
\[
\underline{u} = v + w,
\]
it clearly relies on finding a map $w\in {\cal C}^{2,\alpha}(B^m(a,r),\Bbb{R}^n)$ such that
\begin{equation}\label{3.3}
\left\{ \begin{array}{cccl}
-\Delta w & = & \Gamma(v+w)(\nabla (v+w)\otimes \nabla (v+w)) & \hbox{on }B^m(a,r)\\
w & = & 0 & \hbox{on }\partial B^m(a,r).
\end{array}\right.
\end{equation}
Let us denote by ${\cal C}^{(-\beta)}_{k,\alpha,r}:= {\cal C}^{(-\beta)}_{k,\alpha;B^m(a,r)}$.
We will construct $w$ in ${\cal C}^{(-\alpha)}_{2,\alpha,r}$ by using a fixed point argument.
We first observe that Lemma \ref{LemmaGT0} can be rephrased (and specialized by choosing
$\beta = \alpha$) by saying that there exists a continuous
operator, denoted in the following by $(-\Delta)^{-1}$, from ${\cal C}^{(2-\alpha)}_{0,\alpha,r}$
to ${\cal C}^{(-\alpha)}_{2,\alpha,r}$ which to each $f\in {\cal C}^{(2-\alpha)}_{0,\alpha,r}$ associates
the unique solution $\phi \in {\cal C}^{(-\alpha)}_{2,\alpha,r}$ of
\[
\left\{ \begin{array}{cccl}
-\Delta \phi & = & f & \hbox{on }B^m(a,r)\\
\phi & = & 0 & \hbox{on }\partial B^m(a,r).
\end{array}\right.
\]
We will denote by $\delta$ the norm of $(-\Delta)^{-1}$.
Hence $w\in {\cal C}^{(-\alpha)}_{2,\alpha,r}$ is a solution of (\ref{3.3}) if and only if
\[
w = (-\Delta)^{-1}\left( \Gamma(v+w)(\nabla (v+w)\otimes \nabla (v+w))\right).
\]
Note that $v$ does not belong to ${\cal C}^{(-\alpha)}_{2,\alpha,r}$ (because in particular the trace of $v$
on $\partial B^m(a,r)$ does not vanish). However estimates (\ref{3.2}) holds. This leads us to introduce the
set
\[
{\cal E}_{\alpha;r}:= \{f\in {\cal C}^{2, \alpha}(B^m(a,r), \Bbb{R}^n)/|[f]|_{\alpha;r}<\infty\}
\supset {\cal C}^{(-\alpha)}_{2,\alpha,r},
\]
where the inclusion here is a continuous embedding. We then have:
\begin{lemm}\label{Lemma.3.2}
Let $r\in (0,R_1]$ and $w_0$,$w_1,w_2$ and $w_3$ be in ${\cal E}_{\alpha;r}$. Then
$\Gamma(w_0)(\nabla w_2\otimes \nabla w_3)$ and
$\Gamma(w_1)(\nabla w_2\otimes \nabla w_3)\in {\cal C}^{(2-\alpha)}_{0,\alpha;r}$ and
\begin{equation}\label{3.6}
|(\Gamma(w_1)-\Gamma(w_0))(\nabla w_2\otimes \nabla w_3)|^{(2-\alpha)}_{0,\alpha;r}\leq
C_3r^\alpha |[w_1-w_0]|_{\alpha;r}|[w_2]|_{\alpha;r}|[w_3]|_{\alpha;r},
\end{equation}
\begin{equation}\label{3.6bis}
|\Gamma(w_1)(\nabla w_2\otimes \nabla w_3)|^{(2-\alpha)}_{0,\alpha;r}\leq
C_3r^\alpha |[w_1]|_{\alpha;r}|[w_2]|_{\alpha;r}|[w_3]|_{\alpha;r},
\end{equation}
where $C_3$ is a positive constant.
\end{lemm}
{\em Proof of Lemma \ref{Lemma.3.2}} --- It follows from the interpolation inequality (\ref{1.interpol})
that $|D w_a|^{(1-\alpha)}_{0,\alpha;r}\leq \widetilde{C}|[w_a]|_{\alpha;r}$, $\forall a = 0,1,2,3$. Moreover, for
$r\in (0,R_1]\subset (0,\overline{R}]$,
estimate (\ref{1.dgamma}) implies that
\[
|(\Gamma(w_1)-\Gamma(w_0))|^{(0)}_{0,\alpha;r}\leq C_\Gamma |w_1-w_0|^{(0)}_{0,\alpha;r}
\leq \widetilde{C}C_\Gamma|[w_1-w_0]|_{\alpha;r}.
\]
Hence, using also the inequality $|f|^{(\beta+\gamma)}_{0,\alpha;\Omega}\leq |f|^{(\beta)}_{0,\alpha;\Omega}
|g|^{(\gamma)}_{0,\alpha;\Omega}$, $\forall \beta,\gamma\in \Bbb{R}$ such that $\beta+\gamma\geq 0$
(see 6.11 in \cite{GilbargTrudinger}), we obtain that
\[
\begin{array}{ccl}
|(\Gamma(w_1)-\Gamma(w_0))(\nabla w_2\otimes \nabla w_3)|^{(2-2\alpha)}_{0,\alpha;r} & \leq &
|\Gamma(w_1)-\Gamma(w_0)|^{(0)}_{0,\alpha;r}|D w_2|^{(1-\alpha)}_{0,\alpha;r}|D w_3|^{(1-\alpha)}_{0,\alpha;r}\\
& \leq & \widetilde{C}^3C_\Gamma|[w_1-w_0]|_{\alpha;r}|[w_2]|_{\alpha;r}|[w_3]|_{\alpha;r}.
\end{array}
\]
Thus (\ref{3.6}) follows from the preceding inequality and from
\[
|(\Gamma(w_1)-\Gamma(w_0))(\nabla w_2\otimes \nabla w_3)|^{(2-\alpha)}_{0,\alpha;r}\leq
r^\alpha |(\Gamma(w_1)-\Gamma(w_0))(\nabla w_2\otimes \nabla w_3)|^{(2-2\alpha)}_{0,\alpha;r}.
\]
And (\ref{3.6bis}) is a straightforward consequence of (\ref{3.6}) and of (\ref{1.gamma=0}).
\bbox
\noindent
{\em Back to the proof of Theorem \ref{3.ExistenceTheo}} --- Lemma \ref{Lemma.3.2} allows us to
define the operator
\[
\begin{array}{cccc}
T: & {\cal C}^{(-\alpha)}_{2,\alpha,r} & \longrightarrow & {\cal C}^{(-\alpha)}_{2,\alpha,r}\\
& w & \longmapsto & (-\Delta)^{-1}\left( \Gamma(v+w)(\nabla (v+w)\otimes \nabla (v+w))\right)
\end{array}
\]
and (\ref{3.6bis}) implies
\[
|T(w)|^{(-\alpha)}_{2,\alpha;r}\leq \delta C_3 r^\alpha \left( |[v]|_{\alpha;r} + |w|^{(-\alpha)}_{2,\alpha;r}\right)^3.
\]
In particular, letting ${\cal B}_\Lambda:= \{w\in {\cal C}^{(-\alpha)}_{2,\alpha,r}/|w|^{(-\alpha)}_{2,\alpha;r}
\leq \Lambda \}$ (we recall that $\Lambda$ was defined in (\ref{3.2.5})),
we observe that for all $r\in (0,R_1]\cap (0,R_1']$ where $R_1'= (8\delta C_3 \Lambda^2)^{-1/\alpha}$,
i.e.\,such that in particular $\delta C_3 r^\alpha (2\Lambda)^3<\Lambda$,
\[
\forall w\in {\cal B}_{\Lambda},\quad |T(w)|^{(-\alpha)}_{2,\alpha;r}\leq \Lambda,
\]
which means that $T$ maps the closed ball ${\cal B}_{\Lambda}$ into itself.\\

\noindent
Let us now prove that, for $r$ small enough, the restriction of $T$ on ${\cal B}_{\Lambda}$ is also
contracting: writing that, $\forall w,\widetilde{w}\in {\cal B}_{\Lambda}$,
\[
\begin{array}{ccl}
T(w) - T(\widetilde{w}) & = &
(-\Delta)^{-1}\left[ (\Gamma(v+w)-\Gamma(v+\widetilde{w}))(\nabla (v+w)\otimes \nabla (v+w)\right] \\
& & +\ (-\Delta)^{-1}\left(\Gamma(v+\widetilde{w})(\nabla (v+w)\otimes \nabla (w-\widetilde{w})\right) \\
& & +\ (-\Delta)^{-1}\left( \Gamma(v+\widetilde{w})(\nabla (w-\widetilde{w})\otimes \nabla (v+\widetilde{w})\right) 
\end{array}
\]
and using (\ref{3.6}) we obtain, assuming that $r\leq R_1$,
\[
\begin{array}{ccl}
|T(w) -T(\widetilde{w})|^{(-\alpha)}_{2,\alpha;r} & \displaystyle 
\leq & \delta C_3 r^\alpha |w-\widetilde{w}|^{(-\alpha)}_{2,\alpha;r}\left( |v+w|^{(-\alpha)}_{2,\alpha;r}\right)^2\\
& & +\ \delta C_3 r^\alpha |w-\widetilde{w}|^{(-\alpha)}_{2,\alpha;r}
|v+w|^{(-\alpha)}_{2,\alpha;r}|v+\widetilde{w}|^{(-\alpha)}_{2,\alpha;r}\\
& & +\ \delta C_3 r^\alpha |w-\widetilde{w}|^{(-\alpha)}_{2,\alpha;r}
\left( |v+\widetilde{w}|^{(-\alpha)}_{2,\alpha;r}\right)^2\\
& \displaystyle \leq & 3\delta C_3 r^\alpha \left( 2\Lambda\right)^2|w-\widetilde{w}|^{(-\alpha)}_{2,\alpha;r}.
\end{array}
\]
Hence $T$ is contracting if we further assume that $r<R_1''$, where $R_1'':= (12\delta C_3 \Lambda^2)^{-1/\alpha}$,
because it implies that $3\delta C_3 r^\alpha \left( 2\Lambda\right)^2 < 1$. In conclusion
(observing that actually $R_1''<R_1'$)
if we let $R_2:= \inf \left( R_1, R_1''\right)$,
then for all $r\in (0,R_2) $, $T$ maps the closed ball ${\cal B}_{\Lambda}$ into itself and is contracting.
Hence it admits a unique fixed point $w\in {\cal B}_{\Lambda}$ which is a solution of (\ref{3.3}). \bbox

\section{A maximum principle}
\begin{theo}\label{Theo.4}
Let $({\cal N},h)$ be a pseudo-Riemannian manifold of class ${\cal C}^2$ and $M_0$ be a point in 
${\cal N}$. There exists an open neighbourhood $U_{M_0}$ of $M_0$, a local chart $\phi:U_{M_0}\longrightarrow \Bbb{R}^n$
and constant $\alpha>0$ such that
for any open subset $\Omega$ of $\Bbb{R}^m$ and for any pair of harmonic mappings $u,v:\Omega\longrightarrow (U_{M_0},h_{ij})$
(i.e.\,which satisfy (\ref{0.EquaFonda})), then the function $f:\Omega\longrightarrow \Bbb{R}$ defined 
(using $u\simeq \phi\circ u$ and $v\simeq \phi\circ v$) by
\begin{equation}\label{4.1}
f(x):= (\alpha^2 + |u(x)|^2)(\alpha^2 + |v(x)|^2) {|u(x)-v(x)|^2\over 2},\quad \forall x\in \Omega,
\end{equation}
satisfies the inequality
\begin{equation}\label{4.2}
-\hbox{div} (\rho \nabla f)\leq 0, \quad \hbox{on }\Omega,
\end{equation}
where
\[
\rho(x):= {1\over \alpha^2 + |u(x)|^2}{1\over \alpha^2 + |v(x)|^2},\quad \forall x\in \Omega.
\]
\end{theo}
{\em Remark} ---
Note that here $|\cdot |$ is an Euclidean norm on $U_{M_0}$ which has nothing to do
with the metric $h$ on ${\cal N}$. More precisely, assuming that $\phi(M_0)=0$, for any points $M,\widetilde{M}\in U_{M_0}$
we set $\langle M,\widetilde{M}\rangle := \langle \phi(M)-\phi(M_0), \phi(\widetilde{M})-\phi(M_0)\rangle
= \langle \phi(M), \phi(\widetilde{M})\rangle$, $|M|^2:= |\phi(M)|^2=\langle \phi(M),\phi(M)\rangle$ and
$|M-\widetilde{M}|^2:=|\phi(M)-\phi(\widetilde{M})|^2$.\\

\noindent
{\em Proof of Theorem \ref{Theo.4}} --- Again we first apply Lemma \ref{1.LemmaPrelim} around $M_0$:
it provides us with a local chart $\phi:U_{M_0}'\longrightarrow \Bbb{R}^n$ such that $\phi(M_0)=0$ and
estimates (\ref{1.gamma=0}) and (\ref{1.dgamma}) on the Christoffel symbols $\Gamma^i_{jk}$ hold.
We fix some $\alpha\in (0,\infty)$ which is temporarily arbitrary and whose value will be chosen later.
Then given a pair of harmonic maps $u,v:\Omega\longrightarrow (U_{M_0}',h_{ij})$ we compute $\hbox{div} (\rho \nabla f)$,
where $f$ is given by (\ref{4.1}). We first find that
\[
\rho\nabla f = \langle u-v,\nabla(u-v)\rangle + |u-v|^2\left(
{\langle u,\nabla u\rangle \over \alpha^2+|u|^2} + {\langle v,\nabla v\rangle \over \alpha^2+|v|^2}\right).
\]
Hence (by using the notations $\langle \cdot ,\cdot \rangle$ for the scalar product in $\Bbb{R}^n$ and
$\cdot$ for the scalar product in $\Bbb{R}^m$)
\[
\begin{array}{ccl}
\hbox{div} (\rho \nabla f) & = & |\nabla (u-v)|^2 + \langle u-v,\Delta(u-v)\rangle \\
& & \displaystyle +\ 2\langle u-v,\nabla(u-v)\rangle \cdot \left(
{\langle u,\nabla u\rangle \over \alpha^2+|u|^2} + {\langle v,\nabla v\rangle \over \alpha^2+|v|^2}\right)\\
& & \displaystyle +\  |u-v|^2\left( {|\nabla u|^2\over \alpha^2+|u|^2} + {|\nabla v|^2\over \alpha^2+|v|^2}
+ {\langle u,\Delta u\rangle \over \alpha^2+|u|^2} + {\langle v,\Delta v\rangle \over \alpha^2+|v|^2} \right. \\
& & \displaystyle \ \left. \quad \quad \quad -2 {|\langle u,\nabla u\rangle |^2\over (\alpha^2+|u|^2)^2}
-2 {|\langle v,\nabla v\rangle |^2\over (\alpha^2+|v|^2)^2}\right)\\
& = & G_1 + G_2 +  B_1 + B_2 + B_3 + B_4,
\end{array}
\]
where the ``good'' terms are
\[
G_1:= |\nabla (u-v)|^2,\quad G_2:=
|u-v|^2\left( {|\nabla u|^2\over \alpha^2+|u|^2} + {|\nabla v|^2\over \alpha^2+|v|^2}\right),
\]
and the ``bad'' terms are
\[
\begin{array}{ccl}
B_1 & := & \langle u-v,\Delta(u-v)\rangle \\
B_2 & := & \displaystyle 2\langle u-v,\nabla(u-v)\rangle \cdot \left(
{\langle u,\nabla u\rangle \over \alpha^2+|u|^2} + {\langle v,\nabla v\rangle \over \alpha^2+|v|^2}\right)\\
B_3 & := & \displaystyle |u-v|^2\left( {\langle u,\Delta u\rangle \over \alpha^2+|u|^2} + {\langle v,\Delta v\rangle \over \alpha^2+|v|^2}\right)\\
B_4 & := & \displaystyle -2 |u-v|^2\left( {|\langle u,\nabla u\rangle |^2\over (\alpha^2+|u|^2)^2}
+ {|\langle v,\nabla v\rangle |^2\over (\alpha^2+|v|^2)^2}\right) .
\end{array}
\]
We now need to estimate the bad terms in terms of the good ones.
We let $R\in (0,\infty)$ such that $\phi(U_{M_0}')\subset B^n(0,R)$.
We shall assume in the following that
\begin{equation}\label{4.3}
|u|\leq r\quad \hbox{and}\quad |v|\leq r\quad \hbox{for some }r\in(0,R),
\end{equation}
where $r$ has not yet been fixed. In the following we will first choose $\alpha$ in function of $C_\Gamma$
and $R$, and second we will choose $r$ in function of $\alpha$, $C_\Gamma$ and $R$.\\

\noindent
{\bf Estimation of $B_1$}\\
We have
\[
\begin{array}{ccl}
- \Delta(u-v) & = & \Gamma(u)(\nabla u\otimes \nabla u) - \Gamma (v)(\nabla v\otimes \nabla v)\\
& = & \displaystyle \Gamma(u)(\nabla u\otimes \nabla (u-v)) + \Gamma(u)(\nabla (u-v)\otimes \nabla v)\\
&& + \ (\Gamma(u) - \Gamma(v))(\nabla v\otimes \nabla v).
\end{array}
\]
And because of (\ref{1.gamma=0}) and (\ref{1.dgamma}) which implies $|\Gamma(y)|\leq C_\Gamma|y|$
and $|\Gamma(y) - \Gamma(y')|\leq C_\Gamma|y-y'|$ on $U_{M_0}'$, we deduce that
\[
\begin{array}{ccl}
|\Delta(u-v)| & \leq & C_\Gamma |u|\left( |\nabla u|\,|\nabla(u-v)| + |\nabla v|\, |\nabla(u-v)|\right)\\
&& + \ C_\Gamma |u-v|\, |\nabla v|^2.
\end{array}
\]
Using (\ref{4.3}) and a symmetrisation in $u$ and $v$ one is led to
\[
|\Delta(u-v)|\leq C_\Gamma R|\nabla(u-v)| \left( |\nabla u|+ |\nabla v|\right) +
{C_\Gamma\over 2}|u-v|\left( |\nabla u|^2+ |\nabla v|^2\right).
\]
Hence we deduce using Young's inequality that
\[
\begin{array}{ccl}
|B_1| & \leq & \displaystyle C_\Gamma R|u-v|\,|\nabla(u-v)| \left( |\nabla u|+ |\nabla v|\right)
+ {C_\Gamma\over 2}|u-v|^2\left( |\nabla u|^2+ |\nabla v|^2\right)\\
& \leq & \displaystyle {|\nabla (u-v)|^2\over 4} + 2C^2_\Gamma R^2|u-v|^2\left( |\nabla u|^2+ |\nabla v|^2\right)\\
& & \displaystyle +\ {C_\Gamma\over 2}|u-v|^2\left( |\nabla u|^2+ |\nabla v|^2\right).
\end{array}
\]
We choose $\alpha\in (0,\infty)$ sufficiently small so that ${1\over 2\alpha^2}\geq 4\left(2C^2_\Gamma R^2+{C_\Gamma\over 2}\right)$
and we impose also that $r\leq \alpha$. Then by (\ref{4.3})
\[
|u|,\, |v|\leq r\leq \alpha\quad \Longrightarrow \quad
{1\over \alpha^2+|u|^2}, {1\over \alpha^2+|v|^2}\geq {1\over 2\alpha^2}\geq 4\left(2C^2_\Gamma R^2+{C_\Gamma\over 2}\right)
\]
and thus
\begin{equation}\label{4.4}
|B_1|\leq {G_1\over 4} + {|u-v|^2\over 4}\left( {|\nabla u|^2\over \alpha^2+|u|^2}
+ {|\nabla v|^2\over \alpha^2+|v|^2}\right) = {G_1\over 4} + {G_2\over 4}.
\end{equation}
{\bf Estimation of $B_2$}\\
Using again Young's inequality we obtain
\[
\begin{array}{ccl}
|B_2| & \leq & \displaystyle 2|u-v|\,|\nabla (u-v)|\,\left( {|u|\,|\nabla u|\over \alpha^2+|u|^2}
+ {|v|\,|\nabla v|\over \alpha^2+|v|^2}\right) \\
& \leq & \displaystyle {|\nabla (u-v)|^2\over 2} + 2|u-v|^2\left({|u|\,|\nabla u|\over \alpha^2+|u|^2}
+ {|v|\,|\nabla v|\over \alpha^2+|v|^2}\right)^2 \\
& \leq & \displaystyle {G_1\over 2} + |u-v|^2\left({4|u|^2|\nabla u|^2\over (\alpha^2+|u|^2)^2}
+ {4|v|^2|\nabla v|^2\over (\alpha^2+|v|^2)^2}\right).
\end{array}
\]
We further impose that $r\leq {\alpha\over 4}$. Then by (\ref{4.3})
\begin{equation}\label{4.45}
{4|u|^2\over \alpha^2+|u|^2}, {4|v|^2\over \alpha^2+|v|^2}\leq {4r^2\over \alpha^2}\leq {1\over 4}
\end{equation}
and
\begin{equation}\label{4.5}
|B_2|\leq {G_1\over 2} + |u-v|^2\left({1\over4}{|\nabla u|^2\over \alpha^2+|u|^2}
+ {1\over4}{|\nabla v|^2\over \alpha^2+|v|^2}\right) = {G_1\over 2} + {G_2\over 4}.
\end{equation}
{\bf Estimation of $B_4$}\\
We first write
\[
|B_4|\leq 2|u-v|^2\left({|u|^2\over \alpha^2+|u|^2}{|\nabla u|^2\over \alpha^2+|u|^2}
+ {|v|^2\over \alpha^2+|v|^2}{|\nabla v|^2\over \alpha^2+|v|^2}\right),
\]
and using the fact that ${|u|^2\over \alpha^2+|u|^2}$, ${|v|^2\over \alpha^2+|v|^2}\leq {1\over 16}$ because
of (\ref{4.45}), we deduce that
\begin{equation}\label{4.6}
|B_4|\leq {|u-v|^2\over 8}\left({|\nabla u|^2\over \alpha^2+|u|^2} + {|\nabla v|^2\over \alpha^2+|v|^2}\right)
= {G_2\over 8}.
\end{equation}
{\bf Estimation of $B_3$}\\
We use that
\[
|\Delta u| = |\Gamma(u)(\nabla u\otimes\nabla u)|\leq C_\Gamma|u|\,|\nabla u|^2
\]
and thus $|\langle u,\Delta u\rangle|\leq |u|\,|\Delta u|\leq C_\Gamma|u|^2|\nabla u|^2$
and similarly $|\langle v,\Delta v\rangle|\leq C_\Gamma|v|^2|\nabla v|^2$. Hence by (\ref{4.3})
\[
\begin{array}{ccl}
|B_3| & \leq & \displaystyle C_\Gamma|u-v|^2\left(|u|^2{|\nabla u|^2\over \alpha^2+|u|^2}
+ |v|^2{|\nabla v|^2\over \alpha^2+|v|^2}\right)\\
& \leq & \displaystyle C_\Gamma r^2|u-v|^2\left({|\nabla u|^2\over \alpha^2+|u|^2}
+ {|\nabla v|^2\over \alpha^2+|v|^2}\right).
\end{array}
\]
We further require on $r$ that $C_\Gamma r^2\leq {1\over 4}$. Then
\begin{equation}\label{4.7}
|B_3|\leq {G_2\over 4}.
\end{equation}
{\bf Conclusion}\\
By choosing
\begin{equation}\label{4.75}
\alpha\leq {1\over 2\sqrt{2C_\Gamma^2R^2+C_\Gamma}},\quad
r\leq \inf \left(R,{\alpha\over 4}, {1\over 2\sqrt{C_\Gamma}}\right),
\end{equation}
we obtain using (\ref{4.4}), (\ref{4.5}), (\ref{4.6}) and (\ref{4.7}) that 
\[
\begin{array}{ccl}
\hbox{div} (\rho \nabla f) & = & G_1 + G_2 +  B_1 + B_2 + B_3 + B_4\\
& \geq & \displaystyle G_1 + G_2 - \left({G_1\over 4} + {G_2\over 4}\right) - \left({G_1\over 2} + {G_2\over 4}\right)
- {G_2\over 4} - {G_2\over 8}\\
& = & \displaystyle {G_1\over 4} + {G_2\over 8} \geq 0.
\end{array}
\]
Hence (\ref{4.2}) follows by choosing $U_{M_0}:= \{M\in U_{M_0}'/|\phi(M)|<r\}$ where
$r$ satisfies (\ref{4.75}). \bbox

\section{A result related to capacity}
We prove here the following result.
\begin{lemm}\label{5.Lemma}
Let $\Omega$ be an open subset of $\Bbb{R}^m$, for $m\geq 2$. Let $\rho\in {\cal C}^1(\Omega,\Bbb{R})$ be a function
satisfying $0<A\leq \rho \leq B<\infty$. Let $a\in \Omega$, $\varepsilon_0>0$ such that
$\overline{B^m(a,\varepsilon_0)}\subset \Omega$ and, for all $\varepsilon\in (0,\varepsilon_0)$,
$\Omega_\varepsilon := \Omega\setminus \overline{B^m(a,\varepsilon)}$.\\
Let $\left(\phi_\varepsilon \right)_{\varepsilon\in (0,\varepsilon_0)}$ be a family of functions
$\phi_\varepsilon\in {\cal C}^2(\Omega_\varepsilon )\cap {\cal C}^0(\overline{\Omega_\varepsilon})$ such that
\begin{equation}\label{5.1}
\left\{ \begin{array}{cccl}
\phi_\varepsilon & = & M & \hbox{on }\partial B^m(a,\varepsilon)\\
\phi_\varepsilon & = & 0 & \hbox{on }\partial \Omega\\
-\hbox{div}\left( \rho\nabla \phi_\varepsilon \right) & = & 0 & \hbox{on }\Omega_\varepsilon,
\end{array}\right.
\end{equation}
where $M>0$ is a constant independant of $\varepsilon$. Then for all compact $K\subset \Omega\setminus \{a\}$
and for $\varepsilon\in (0,\varepsilon_0)$ such that $K\subset \Omega_\varepsilon$, the restriction of $\phi_\varepsilon$
on $K$ converges to 0 in $L^1(K)$ when $\varepsilon$ tends to 0.
\end{lemm}
{\em Proof} --- A first step consists in proving that the energy of $\phi_\varepsilon$,
\[
E_\varepsilon:= {\cal A}_\varepsilon[\phi_\varepsilon]:= 
\int_{\Omega_\varepsilon} \rho|\nabla \phi_\varepsilon|^2\, dx
\]
converges to 0 when $\varepsilon$ tends to 0. This is a very standard result which can be checked as follows:
we know that $\phi_\varepsilon$ is energy minimizing and hence that
$E_\varepsilon\leq {\cal A}_\varepsilon[f]$, for all $f\in {\cal C}^2(\Omega_\varepsilon )\cap {\cal C}^0(\overline{\Omega_\varepsilon})$
such that $f = M$ on $\partial B^m(a,\varepsilon)$ and $f = 0$ on $\partial \Omega$.
One can choose for $f$ $v_\varepsilon(x) := M\chi(x)G_\varepsilon(x-a)$,
where $\chi \in {\cal C}^2(\Omega_\varepsilon)$ satisfies $0\leq \chi\leq 1$,
$\hbox{supp}\chi\subset B^m(a,\varepsilon_0)$, $\chi=1$ on $B^m(a,\varepsilon_0/2)$ and $|\nabla \chi|\leq 4/\varepsilon_0$
and where $G_\varepsilon$ is the Green function on $\Bbb{R}^m$
($G_\varepsilon(x):= {\log |x|\over \log \varepsilon}$ if $m=2$ and $G_\varepsilon(x):= {\varepsilon^{m-2}\over |x|^{m-2}}$
if $m\geq 3$). Then a straightforward computation shows that $\lim_{\varepsilon\rightarrow 0}{\cal A}_\varepsilon[v_\varepsilon]=0$.
Hence
\begin{equation}\label{5.2}
\lim_{\varepsilon\rightarrow 0}E_\varepsilon = 0.
\end{equation}
Second we consider, for all $s\in [0,M]$, the level sets
\[
\Omega_\varepsilon^s:= \{x\in \Omega_\varepsilon/ \phi_\varepsilon(x)>s\}.
\]
Note that the maximum principle implies that $\phi_\varepsilon$ takes values in $[0,M]$. Sard's Theorem
implies that the set $V_c:= \{s\in [0,M]/\exists x\in \Omega_\varepsilon \hbox{ such that }\phi_\varepsilon(x)=s
\hbox{ and } \nabla \phi_\varepsilon(x)= 0\}$ (critical values) is negligeable (moreover it is also closed\footnote{and
$V_c\subset (0,M)$ because of the Hopf maximum principle}). And
$\forall s\in [0,M]\setminus V_c$, $\partial \Omega_\varepsilon^s=\{ x\in \Omega_\varepsilon/ \phi_\varepsilon(x)=s\}$
is a smooth submanifold. We let $\Gamma^s$ be the exterior part of $\partial \Omega_\varepsilon^s$ so that
we have the splitting
\[
\partial \Omega_\varepsilon^s = \Gamma^s\cup \partial B^m(a,\varepsilon).
\]
Using the equation (\ref{5.1}) we observe that, $\forall s,s'\in [0,M]\setminus V_c$ such that $s<s'$,
\[
\begin{array}{ccl}
0 & = & \displaystyle \int_{\Omega_\varepsilon^s\setminus \Omega_\varepsilon^{s'}}
- \hbox{div}\left(\rho\nabla \phi_\varepsilon\right) dx\\
 & = & \displaystyle \int_{\Gamma^s}-\rho\langle \nabla \phi_\varepsilon,\nu\rangle d{\cal H}^{m-1}
 + \int_{\Gamma^{s'}}\rho\langle \nabla \phi_\varepsilon,\nu\rangle d{\cal H}^{m-1}\\
 & = & \displaystyle \int_{\Gamma^s}\rho|\nabla \phi_\varepsilon|d{\cal H}^{m-1}
 - \int_{\Gamma^{s'}}\rho|\nabla \phi_\varepsilon|d{\cal H}^{m-1},
\end{array}
\]
where $d{\cal H}^{m-1}$ is the $(m-1)$-dimensional Hausdorff measure. Here
we have used in the last line the fact that, on $\Gamma^s$ and $\Gamma^{s'}$,
$\nabla \phi_\varepsilon$ is parallel (and of opposite orientation) to the normal vector $\nu$.
This implies that the function
\begin{equation}\label{5.3}
[0,M]\setminus V_c\ni s\longmapsto \int_{\Gamma^s}\rho|\nabla \phi_\varepsilon|d{\cal H}^{m-1}
\quad \hbox{is constant.}
\end{equation}
We now use the coarea formula to obtain
\[
\begin{array}{ccl}
\displaystyle \int_{\Omega_\varepsilon}\rho|\nabla \phi_\varepsilon|^2 dx& = & 
\displaystyle \int_0^Mds\int_{\Gamma^s}\rho|\nabla \phi_\varepsilon|^2{d{\cal H}^{m-1}\over |\nabla \phi_\varepsilon|}\\
& = & \displaystyle \int_0^Mds\int_{\Gamma^s}\rho|\nabla \phi_\varepsilon|d{\cal H}^{m-1}.
\end{array}
\]
Thus we deduce that, using (\ref{5.3}),
\begin{equation}\label{5.4}
\forall s\in [0,M]\setminus V_c,\quad \int_{\Gamma^s}\rho|\phi_\varepsilon|d{\cal H}^{m-1}
= {1\over M}\int_{\Omega_\varepsilon}\rho|\nabla \phi_\varepsilon|^2 dx = {E_\varepsilon\over M}.
\end{equation}
We now let $F_\varepsilon:[0,M]\longrightarrow [0,\infty)$ be the function defined by
\[
F_\varepsilon(s):= |\Omega_\varepsilon^s\cup \overline{B^m(a,\varepsilon)}|, 
\quad \hbox{the Lebesgue measure of }\Omega_\varepsilon^s\cup \overline{B^m(a,\varepsilon)}.
\]
Obviously $F_\varepsilon$ is a decreasing function and so $F_\varepsilon'$ is a nonpositive measure.
We can decompose this measure as $F'_\varepsilon = (F'_\varepsilon)_a + (F'_\varepsilon)_s$,
where $(F_\varepsilon')_a$ is the absolutely continuous part of $F_\varepsilon'$ and
$(F'_\varepsilon)_s$ is the singular part of $F_\varepsilon'$. Moreover $F_\varepsilon$ is differentiable
on $[0,M]\setminus V_c$ with
\[
\forall s\in [0,M]\setminus V_c,\quad F_\varepsilon'(s) = - 
\int_{\Gamma^s}{d{\cal H}^{m-1}\over |\nabla \phi_\varepsilon|}
\]
and $\hbox{supp}(F_\varepsilon')_s\subset V_c$. We deduce from this identity and
from (\ref{5.4}), by using the Cauchy--Schwarz inequality, that $\forall s\in [0,M]\setminus V_c$,
\[
\begin{array}{ccl}
|\Gamma^s| & = & \displaystyle \int_{\Gamma^s}d{\cal H}^{m-1} = 
\int_{\Gamma^s}\sqrt{|\nabla \phi_\varepsilon|}{d{\cal H}^{m-1}\over \sqrt{|\nabla \phi_\varepsilon|}}\\
& \leq & \displaystyle \sqrt{\int_{\Gamma^s}|\nabla \phi_\varepsilon|d{\cal H}^{m-1}}
\sqrt{\int_{\Gamma^s}{d{\cal H}^{m-1}\over |\nabla \phi_\varepsilon|}}\\
& \leq & \displaystyle \sqrt{{1\over A}\int_{\Gamma^s}\rho|\nabla \phi_\varepsilon|d{\cal H}^{m-1}}
\sqrt{-F'_\varepsilon(s)}\\
& \leq & \displaystyle \sqrt{E_\varepsilon\over AM}\sqrt{-F'_\varepsilon(s)}.
\end{array}
\]
Hence
\begin{equation}\label{5.5}
\forall s\in [0,M]\setminus V_c,\quad |\Gamma^s|^2 \leq -{E_\varepsilon\over AM}F'_\varepsilon(s).
\end{equation}
We observe that this inequality extends on the whole interval $[0,M]$ in the sense of measure: $V_c$
is Lebesgue negligeable and if $s$ is a singular point of $F'_\varepsilon$, then the above
inequality holds since the left hand side is a function.
We next exploit (\ref{5.5}) together with the isoperimetric inequality for the
subset $\Omega_\varepsilon^s\cup \overline{B^m(a,\varepsilon)}\subset \Bbb{R}^m$ and its
boundary $\Gamma^s$:
\begin{equation}\label{5.6}
m^{m-1}\omega_mF_\varepsilon(s)^{m-1} =
m^{m-1}\omega_m|\Omega_\varepsilon^s\cup \overline{B^m(a,\varepsilon)}|^{m-1} \leq |\Gamma^s|^m.
\end{equation}
Then (\ref{5.5}) and (\ref{5.6}) imply
\begin{equation}\label{5.7}
F'_\varepsilon + k_\varepsilon F_\varepsilon^{2(m-1)/m}\leq 0,
\quad \hbox{with }k_\varepsilon:= {AM\left(m^{m-1}\omega_m\right)^{2/m}\over E_\varepsilon},
\end{equation}
in the sense of measure on $[0,M]$.\\
\noindent
{\bf The case $m=2$}\\
Equation (\ref{5.7}) then implies
\[
\forall s\in [0,M],\quad F_\varepsilon(s)\leq F_\varepsilon(0)e^{-k_\varepsilon s} = |\Omega|e^{-k_\varepsilon s}.
\]
Hence for any compact $K\subset \Omega\setminus \{a\}$ and for $\varepsilon$ small enough, by using
the coarea formula, we have
\[
||\phi_\varepsilon||_{L^1(K)} = \int_0^M|K\cap \Omega_\varepsilon^s|ds
\leq \int_0^MF_\varepsilon(s)ds\leq {|\Omega|\over k_\varepsilon}\left(1-e^{-k_\varepsilon M}\right).
\]
which implies that $||\phi_\varepsilon||_{L^1(K)}$ tends to 0 when $\varepsilon\rightarrow 0$
because $k_\varepsilon$ tends to $\infty$, because of (\ref{5.2}). Hence the Lemma is
proved in this case.\\

\noindent
{\bf The case $m\geq 2$}\\
Let us denote by $\beta:= 2{m-1\over m}-1\in (0,1)$. We deduce analogously to the preceding case that
\[
\forall s\in [0,M],\quad F_\varepsilon(s)\leq {|\Omega|\over (1+\beta|\Omega|^\beta k_\varepsilon s)^{1/\beta}}
\]
and thus
\[
||\phi_\varepsilon||_{L^1(K)}\leq {|\Omega|^{1-\beta}\over (1-\beta)k_\varepsilon}
\left(1-{1\over \left(1+\beta|\Omega|^\beta k_\varepsilon M\right)^{1/\beta-1}}\right),
\]
which leads to the same conclusion. \bbox

\section{The proof of the Main theorem}
We conclude this paper by proving Theorem \ref{0.MainTheo}. Let $u$ be a continuous map
from $\Omega$ to ${\cal N}$ and assume that $u$ is ${\cal C}^2$ and harmonic with values in
$({\cal N},h_{ij})$ on $\Omega\setminus \{a\}$.
Using Theorem \ref{Theo.4} with $M_0=u(a)$ we deduce that there exists a
neighbourhood $U_{u(a)}$ of $u(a)$ in ${\cal N}$ such (\ref{4.2}) holds for any
pair of harmonic maps into $(U_{\phi(a)},h_{ij})$. We hence can restrict $u$ to a ball
$B^m(a,R)$, where $R$ is chosen so that $u(B^m(a,R))\subset U_{u(a)}$. Then we use
the existence result \ref{3.ExistenceTheo}: we deduce that
there exists some $R_2\in (0,R)$ such that, for any $r\in (0,R_2)$ there exists a 
map $\underline{u}\in {\cal C}^{2,\alpha}(B^m(a,r))\cap {\cal C}^{0,\alpha}(\overline{B^m(a,r)})$
which is harmonic into $(U_{u(a)},h_{ij})$ and which coincides with $u$ on $\partial B^m(a,r)$.
Then we choose some $r\in (0,R_2)$ and we identify
$u\simeq \phi\circ u$ and $\underline{u}\simeq \phi\circ \underline{u}$ as in Theorem \ref{Theo.4}.
Note that it is clear that there exists some $A\in (0,\infty)$ such that $|u|,|\underline{u}|\leq A$
on $B^m(a,r)$. Now let
\[
f(x):= (\alpha^2 + |u(x)|^2)(\alpha^2 + |\underline{u}(x)|^2) {|u(x)-\underline{u}(x)|^2\over 2},
\quad \forall x\in B^m(a,r).
\]
where $\alpha$ has been chosen as in Theorem \ref{Theo.4}. For any
$\varepsilon>0$ such that $\overline{B^m(a,\varepsilon)}\subset B^m(a,r)$, we consider the
map $\phi_\varepsilon\in {\cal C}^{2}(B^m(a,r)\setminus \overline{B^m(a,\varepsilon)})\cap
{\cal C}^{0}(\overline{B^m(a,r)}\setminus B^m(a,\varepsilon))$ which is the
solution to
\[
\left\{ \begin{array}{cccl}
\phi_\varepsilon & = & M & \hbox{on }\partial B^m(a,\varepsilon)\\
\phi_\varepsilon & = & 0 & \hbox{on }\partial B^m(a,r)\\
-\hbox{div}\left( \rho\nabla \phi_\varepsilon \right) & = & 0 & \hbox{on }B^m(a,r)\setminus \overline{B^m(a,\varepsilon)},
\end{array}\right.
\]
where
\[
\rho(x):=  {1\over \alpha^2 + |u(x)|^2}{1\over \alpha^2 + |\underline{u}(x)|^2},\quad \forall x\in B^m(a,r)
\]
and
\[
M:= 2A^2(\alpha^2 + A^2)(\alpha^2 + A^2).
\]
Clearly we have $f\leq \phi_\varepsilon$ on $\partial (B^m(a,r)\setminus \overline{B^m(a,\varepsilon)})$
and Theorem \ref{Theo.4} implies that
$-\hbox{div}\left( \rho\nabla f \right) \leq 0 = -\hbox{div}\left( \rho\nabla \phi_\varepsilon \right)$
on $B^m(a,r)\setminus \overline{B^m(a,\varepsilon)}$. Hence the maximum principle implies that
$f\leq \phi_\varepsilon$ on $B^m(a,r)\setminus \overline{B^m(a,\varepsilon)}$. Now if we fix a compact
subset $K\subset B^m(a,r)\setminus \{a\}$ and suppose that $\varepsilon$ is sufficiently small so that $K\subset B^m(a,r)\setminus 
\overline{B^m(a,\varepsilon)}$, the inequality $f\leq \phi_\varepsilon$ on $K$ implies
\[
||f||_{L^1(K)}\leq ||\phi_\varepsilon||_{L^1(K)}.
\]
Letting $\varepsilon$ tend to 0 and using Lemma \ref{5.Lemma} we deduce that $||f||_{L^1(K)}=0$.
Since $K$ is arbitrary and $f$ is continuous on $B^m(a,r)$, we conclude $f=0$ on $B^m(a,r)$. Hence
$u$ coincides with $\underline{u}$ on $B^m(a,r)$. Thus $u$ is ${\cal C}^{2,\alpha}$ on $B^m(a,r)$.
\bbox

\noindent
{\em Fr\'ed\'eric H\'elein\\
CMLA, ENS de Cachan, 61 avenue du Pr\'esident Wilson,\\
94235 Cachan cedex, France\\
helein@cmla.ens-cachan.fr}


\begin{thebibliography}{100}

\bibitem{Bryant} R. Bryant, {\em A duality theorem for Willmore surfaces}, J. Differential Geom. 20
(1984), 23--53.

\bibitem{BFLPP} F. Burstall, D. Ferus, K. Leschke, F. Pedit, U. Pinkall, {\em Conformal geometry of
surfaces in $S^4$ and quaternions}, Lecture Notes in Mathematics, Springer, Berlin, Heidelberg, 2002.

\bibitem{GilbargTrudinger} D. Gilbarg, N.S. Trudinger, {\em Elliptic partial differential equations
of second order}, Springer-Verlag Berlin Heidelberg 2001.

\bibitem{Helein} F. H\'elein, {\em Regularity and uniqueness of harmonic maps into an ellipsoid},
Manuscripta Math. 60 (1988), 235--257.

\bibitem{Helein2} F. H\'elein, {\em Willmore immersions and loop groups}, J. Differential Geom. 50
(1998), 331--385.

\bibitem{JaegerKaul} W. J\"ager, H. Kaul, {\em Uniqueness and stability of harmonic maps and
their Jacobi fields}, Manuscripta Math. 28 (1979), 269--291.

\bibitem{LU} O. Ladyzhenskaya, N. Ural'tseva, {\em Linear and quasilinear elliptic equations}, Academic
Press, 1968.

\bibitem{Morrey} C.B. Morrey, {\em Multiple integrals in the calculus of variations}, Grundleheren 130,
Springer Berlin, 1966.

\bibitem{SacksUhlenbeck} J. Sacks, K. Uhlenbeck, {\em The existence of minimal immersions of
2-spheres}, Ann. Math. 113 (1981), 1--24.

\bibitem{Stein} E. Stein, {\em Singular integrals and differentiability properties of functions},
Princeton University Press, 1970.


\end{thebibliography}
\end{document}